\documentclass{article}

\usepackage[utf8]{inputenc}
\usepackage[T1]{fontenc}
\usepackage[a4paper]{geometry}
\usepackage[ruled,longend,vlined]{algorithm2e}
\usepackage{amsmath}
\usepackage{amssymb}
\usepackage{hyperref}

\newcommand{\F}{\mathbb {F}}
\newcommand{\fp}{\F_p}
\newcommand{\fq}{\F_q}
\newcommand{\fpk}{\F_{p^k}}
\newcommand{\fqk}{\F_{q^k}}
\newcommand{\fqbar}{\overline\F_q}

\newcommand{\curvefii}{$E_{12, 1}$}
\newcommand{\curvefiv}{$E_{12, 2}$}
\newcommand{\curvefv}{$E_{12, 3}$}
\newcommand{\curvefall}{$E_{12, i}$}
\newcommand{\curveg}{$E_{11}$}
\newcommand{\curveh}{$E_{19}$}
\newcommand{\curvei}{$E_{18}$}
\newcommand{\curvel}{$E_{24}$}
\newcommand{\curvem}{$E_{27}$}
\newcommand{\curven}{$E_{26}$}
\newcommand{\curveo}{$E_{25}$}
\newcommand{\curvep}{$E_{15}$}
\newcommand{\curves}{$E_{16}$}
\newcommand{\curvet}{$E_{17}$}
\newcommand{\curvez}{$E_9$}

\newcommand {\ew}{e_{\mathrm W}}
\newcommand {\et}{e_{\mathrm T}}
\newcommand {\ea}{e_{\mathrm A}}
\newcommand {\etw}{e_{\mathrm {A'}}}
\newcommand {\ev}{e_{\mathrm O}}
\newcommand {\evt}{e_{\mathrm O'}}
\newcommand {\dpol}{F}
\newcommand {\OO}{\mathcal O}
\newcommand {\ord}{\operatorname {ord}}
\newcommand {\ddiv}{\operatorname {div}}
\newcommand {\Ker}{\operatorname {ker}}
\newcommand {\Aut}{\operatorname {Aut}}

\title{Implementing cryptographic pairings \\at standard security levels}
\author{Andreas Enge\footnote {INRIA, LFANT, F-33400 Talence, France \newline
CNRS, IMB, UMR 5251, F-33400 Talence, France \newline
Univ. Bordeaux, IMB, UMR 5251, F-33400 Talence, France \newline
andreas.enge@inria.fr}\, and
Jérôme Milan\footnote{INRIA, CNRS, UMR 7161, École polytechnique, LIX,
91128 Palaiseau, France}}
\date{22 July 2014}

\begin{document}
%-----------------------------------------------------------------------------

\maketitle
\begin{abstract}
This study reports on an implementation of cryptographic pairings
in a general purpose computer algebra system.
For security levels equivalent to the different AES flavours,
we exhibit suitable curves in parametric families and show that
optimal ate and twisted ate pairings exist and can be efficiently
evaluated. We provide a correct description of Miller's algorithm for signed
binary expansions such as the NAF and extend a recent variant due
to Boxall \textit {et al.} to addition-subtraction chains.
We analyse and compare several algorithms
proposed in the literature
for the final exponentiation.
Finally, we give recommendations on which curve and pairing to
choose at each security level.
\end{abstract}

%-----------------------------------------------------------------------------
\section{Pairings on elliptic curves}
\label{section:pairings}
%-----------------------------------------------------------------------------

In this article, we treat cryptographic
bilinear pairings $G_1 \times G_2 \to G_T$ on
elliptic curves~$E$ defined over some finite field~$\fq$ of
characteristic~$p$.
We emphasise that our aim is not to set new speed records for particular
curves, cf. \cite {ArKaLoGeLo11,PeSiNaBa11,ArFuKnMeRo13}, but to
compare different choices of pairings and parameters at various
security levels, using a general purpose, but reasonably optimised,
implementation in a general purpose
computer algebra system. Such an analysis will be meaningful assuming that
the ratios between the various operations remain constant when switching to
hand-optimised assembly implementations in each instance.

We fix the following standard notations and setting.
Let $E (\fq)$ denote the $\fq$-rational points on $E$, and let $r$ be a prime
divisor of $\# E (\fq) = q + 1 - t$ that does not divide $q-1$,
where $t$ is the trace of Frobenius. Let the embedding degree~$k$ be the
smallest integer such that $r$ divides $q^k - 1$, and denote by $\pi$ the
Frobenius map $E (\fqk) \to E (\fqk)$, $(x, y) \mapsto (x^q, y^q)$.
The $r$-torsion subgroup $E [r]$ is defined over $\fqk$, and it contains the
non-trivial subgroup $E (\fq) [r]$ of $\fq$-rational $r$-torsion points.
Denote by $\mu_r$ the subgroup of $r$-th roots of unity in $\fqk^\ast$.

Typically, $G_T = \mu_r$, $G_1 = E (\fq) [r]$, and $G_2$ is a subgroup of
order~$r$ of either $E [r]$ or of $E (\fqk) / r E (\fqk)$.

\subsection {Functions with prescribed divisors}
\label {subsection:divisors}

Let $E$ be given over $\fq$ by an equation
in the variables $x$ and $y$. For a rational function
$f \in \fqbar (E) := \fqbar (x)[y] / (E)$ and a point $P \in E$,
denote by $\ord_P (f)$ the positive multiplicity of the zero $P$ of $f$,
the negative multiplicity of the pole $P$ of $f$, or $0$ if $P$ is neither a
zero nor a pole of $f$. Denote by
$\ddiv (f) = \sum_P \ord_P (f) [P]$ the divisor of $f$, an
element of the free abelian group generated by the symbols~$[P]$, where $P$
is a point on~$E$.

The definition and computation of pairings involve certain rational functions
with given divisors, in particular, $f_{n,P}$ with
\[
\ddiv (f_{n, P}) = n[P]-[nP]-(n-1)[\OO],
\]
the lines $\ell_{P, Q}$ through two (not necessarily distinct) points $P$ and $Q$
with
\[
\ddiv (\ell_{P, Q}) = [P] + [Q] + [-(P+Q)] - 3 [\OO]
\]
and the vertical lines $v_P$ through a point $P$ with
\[
\ddiv (v_P) = [P] + [-P] - 2 [\OO].
\]
All these functions are defined up to a multiplicative
constant, and they are normalised at infinity by the condition
$\left( f \left( \frac {Y}{X} \right)^{\ord_\OO (f)} \right) (\OO) = 1$.

In particular, we have $\ell_{P, -P} = v_P$, $f_{1, P} = 1$ and
$f_{-1, P} = 1 / v_P$.

The function $f_{n, P}$ is of degree $O (n)$ and may be evaluated
in $O (\log n)$ steps by the algorithms of \S\ref {section:loop}.

\subsection {Cryptographic pairings}
\label {subsection:pairings}

We quickly recall the main cryptographic pairings. In applications, they are
usually restricted to $E (\fq)[r]$ in one argument and to a subgroup of
order~$r$ in the other argument.

\paragraph {Weil pairing}
\[
\ew : E(r) \times E(r) \to \mu_r, \quad
      (P, Q) \mapsto (-1)^r \frac{f_{r,P}(Q)}{f_{r,Q}(P)}
\]
Computing the pairing requires the evaluation of two functions; moreover, with
$P \in E(\fq)$ and $Q \in E(\fqk)$,
the function $f_{r,Q}$ is much costlier to evaluate by the algorithms of
\S\ref {section:loop}.

\paragraph{Tate pairing}
\[
\begin {array}{ccccc}
\et : E(\fq)[r] \times E(\fqk)/rE(\fqk) & \to & \fqk^\ast/(\fqk^\ast)^r & \simeq & \mu_r, \\
      (P, Q) & \mapsto & f_{r, P}(Q) & \leftrightarrow & f_{r, P}(Q)^{(q^k-1)/r}.
\end {array}
\]
The pairing requires only one evaluation of a rational function, but
the original definition with a quotient group as domain is unwieldy since
there is no easy way of defining unique representatives. The final
exponentiation step of raising to the power $\frac {q^k-1}{r}$ realises an
isomorphism with $\mu_r$, and the resulting pairing is usually called the
reduced Tate pairing.

\paragraph{Ate pairing}
By restricting its arguments to the two eigenspaces of $E [r]$ under
the Frobenius with eigenvalues $1$ and $q$, respectively,
the ate pairing introduced in~\cite{HeSmVe06} replaces $f_{r, P} (Q)$
with $r \in O (q)$ by $f_{T, Q} (P)$, where
$T = t - 1 \in O (\sqrt q)$. This saving may be offset by the swap of
arguments $P$ and $Q$, so that the function is defined over $\fqk$
instead of $\fq$.
\[
\ea : E(\fq)[r] \times E [r] \cap \Ker(\pi - q) \to \mu_r,
\quad
(P, Q) \mapsto f_{T,Q}(P)^{(q^k-1)/r}.
\]

\paragraph{Twisted ate pairing}
The curve $E$ admits a twist of degree $d=\gcd(k, \#\Aut (E))$;
let $e=k/d$.
The twisted ate pairing of \cite{HeSmVe06} works again with a function
over $\fq$, at the expense of a factor of $O (e)$ in its evaluation:
\[
\etw : E(\fq)[r] \times E [r] \cap \Ker (\pi - q) \to \mu_r,
\quad
(P, Q) \mapsto f_{T^e, P}(Q)^{(q^k-1)/r}.
\]

\paragraph {Optimal pairing}
Generalisations of the ate and twisted ate pairings requiring
several functions for their evaluation have been given in
\cite{Hess08,Vercauteren10}. All of them take $G_1 = E (\fq)$
and $G_2 = E [r] \cap \Ker (\pi - q)$.
They are
evaluated with low degree functions, typically requiring
$O (\log r / \varphi (k))$ operations (not counting the final
exponentiation), where $\varphi$ is Euler's function.

Let $\lambda = mr =\sum_{i=0}^n \lambda_i q^i$ be a suitably chosen
multiple of $r$ with $r \nmid m$ such that the $\lambda_i$ are small;
more precisely, one requires a short addition-subtraction sequence passing
through all $\lambda_i$.
An optimal ate pairing is obtained by
\[
\ev : G_1 \times G_2 \to \mu_r, \quad
(P, Q) \mapsto
\left(
\prod_{i=0}^n   f_{\lambda_i, Q}^{q^i}(P)
\prod_{i=0}^{n-1} \frac{\ell_{s_{i+1} Q, \lambda_i q^i Q}(P)}{v_{s_i Q}(P)}
\right )^{(q^k-1)/r},
\]
where $s_i = \sum_{j=i}^n \lambda_j q^j$.
Since $\lambda = \varphi_k (q)$ yields a degenerate pairing, one may
assume $n < \varphi (k)$; a precise condition for non-degeneracy
is given in \cite[Th.~4]{Vercauteren10}.
Finding such a multiple of $r$ with small coefficients in base $q$
is a common integer
relation problem. It may be solved, for example, by using the LLL algorithm
to find a short vector $(\lambda_0, \ldots, \lambda_{\varphi (k) - 1})$
in the lattice generated by $(r, 0, \ldots, 0)$ and the
$(-q^{i-1}, 0, \ldots, 0, 1, 0, \ldots, 0)$
with $1$ at position $i$ for $i = 2, \ldots, \varphi (k)$.

\paragraph {Optimal twisted pairing}
In the presence of a twist of degree $d$ such that $k = d e$,
a pairing can be obtained in a similar fashion
from $\lambda = mr = \sum_{i=0}^n \lambda_i T^{ei}$
with $n < \varphi (d)$. The only interesting cases are
$d \in \{ 3, 4, 6 \}$ with $\varphi (d) = 2$ (otherwise, we obtain
again the twisted ate pairing). Then
\[
\evt : G_1 \times G_2 \to \mu_r, \quad
(P, Q) \mapsto
\left(
f_{\lambda_0,P}(Q) f_{\lambda_1,P}^{q^e}(Q) v_{\lambda_0 P}(Q)
\right )^{(q^k-1)/r}
\]
defines a pairing, where
$(\lambda_0, \lambda_1)$ is a short vector in the lattice generated by
$(r, 0)$ and $(-T^e, 1)$.

%-----------------------------------------------------------------------------
\section{Curves and associated optimal pairings}
\label {section:curve}
%-----------------------------------------------------------------------------

\subsection {Curve selection}

Pairing-based cryptographic settings should be designed so that the discrete
logarithm problem is intractable in each group involved.
Current best attacks on a subgroup of prime order $r$ of an elliptic curve
require on the order of $\sqrt r$ operations. For the sake of efficiency of
the implementation and to save bandwidth in transmitting points on the curve,
it would be desirable that $\# E (\fq) = r \approx q$ by Hasse's theorem.
This condition, however, is not easily met for arbitrary embedding degree $k$,
and the parameter $\rho = \frac {\log q}{\log r}$ measures to which extent
it is violated: Values close to~$1$ would be optimal.
Discrete logarithms in finite fields~$\fqk$ of large characteristic,
which we deal with in this article, may be computed by algorithms with
subexponential complexity $L_{q^k} (1/3)$
(as opposed to quasi-polynomial complexity in small characteristic
\cite{BaGaJoTh14}).
One is thus looking for parameter values such that
$\sqrt r \approx L_{q^k} (1/3) = L_{r^{\rho k}} (1/3)$.
Taking logarithms on both sides shows that for bounded~$\rho$,
the embedding degree $k$ grows asymptotically as $\log^2 q$.

Several studies have refined this argument.
Table~\ref {table:ecrypt} summarises security levels equivalent to the
AES secret key sizes, all given in bits, according to the Ecrypt2
recommendations~\cite{Smart10}, validated by a large
community of cryptologists.
\begin {table}[!hbt]
\centerline {
\begin{tabular}{c|r|r|r}
Security & $\log_2 r$ & $\log_2 {q^k}$ & Target $k \rho$\\
\hline
128 & 256 & 3248 & 12.7\\
192 & 384 & 7936 & 20.7\\
256 & 512 & 15424& 30.1
\end{tabular}
}
\caption {Recommended curve sizes}
\label{table:ecrypt}
\end {table}

For each security level, we have selected from \cite{FrScTe10}
several curve families approximately fulfilling the requirements on the
parameter sizes,
trying to stay as close as possible to the ideal $k\rho$ values.

These curves, together with their main parameters, are given in
Table~\ref {table:curves}; for more details, see Appendix~\ref{appendix:params}.
The entries in column ``construction'' refer to the algorithms
of~\cite {FrScTe10}.
Supersingular curves are ruled out at these security levels by their
too small embedding degree $k \leq 6$, so we restricted the search to
ordinary curves. Since non-prime base fields $\fq$ are virtually impossible
to reach, all curves are defined over prime fields $\fp$.
We favoured a small Hamming weight of~$r$ and field extensions $\fpk$ which
may be defined by a binomial $X^k-a$ with a small $a \in \fp$.
For comparison purposes, we also included the Barreto-Naehrig curves
\curvefiv\ and \curvefv,
widely considered the best choice for 128 bit security due to their optimal
value of $\rho =1$, for security levels of 192 and 256 bits
by artificially increasing $p$.
Some curve families are very sparse, which forced us to relax the constraints,
for instance for \curveo.
\begin {table}[hbt]
\centerline {
\scriptsize
\begin{tabular}{c|r|r|r|c|c|r|r|r|r}
Security & $k$ & $\rho$  & $k \rho$ & Curve     & Construction & $\deg q$ & $\deg r$ & $\deg t$ & $\varphi(k)$\\
\hline
128      &   9 &  $4/3$  &     12.0 & \curvez   & 6.6          &   8  &   6  &  4   &  6\\
         &  11 &   $6/5$ &    13.2  & \curveg   & 6.6          &  24  &  20  & 12   & 10\\
         &  12 &       1 &    12    & \curvefii & 6.8          &  4   &   4  &  2   &  4\\
\hline
192      &  12 &       1 &    12    & \curvefiv & 6.8          &   4  &   4  &  2   &  4\\
         &  15 &   $3/2$ &    22.5  & \curvep  	& 6.6          &  12  &   8  &  6   &  8\\
         &  16 &   $5/4$ &    20    & \curves   & 6.11         &  10  &   8  &  5   &  8\\
           %&  17 &   $9/8$ &    19.1  & --        & 6.6          &  36  &  32  & 18   & 16\\
         &  17 & $19/16$ &    20.2  & \curvet   & 6.2          &  38  &  32  &  2   & 16\\
         &  18 &   $4/3$ &    24    & \curvei   & 6.12         &   8  &   6  &  4   &  6\\
         &  19 &  $10/9$ &    21.1  & \curveh   & 6.6          &  40  &  36  & 20   & 18\\
\hline
256      &  12 &       1 &    12    & \curvefv  & 6.8          &   4  &   4  &  2   &  4\\
         &  24 &   $5/4$ &    30    & \curvel   & 6.6          &  10  &   8  &  1   &  8\\
         &  25 & $13/10$ &    32.5  & \curveo   & 6.6          &  52  &  40  & 26   & 20\\
         &  26 &   $7/6$ &    30.34 & \curven   & 6.6          &  28  &  24  & 14   & 12\\
         &  27 &  $10/9$ &    30    & \curvem   & 6.6          &  20  &  18  & 10   & 18\\
\end{tabular}
}
\caption {Main curve parameters}
\label{table:curves}
\end {table}

\subsection {Optimal pairings}
\label {ssec:optimal}

Each curve family gives rise to different optimal pairings,
which are obtained via short vectors in certain lattices as
explained in \S\ref{subsection:pairings}.
For each curve family in Table~\ref {table:curves}, such short vectors
are given in Table~\ref {table:optivec}.
Notice that for \curvet\ and \curvel, the ate pairings are already
optimal, and that for $d \not\in \{ 3, 4, 6 \}$ and moreover for
\curvel, the twisted optimal pairings are exactly the twisted ate pairings.

%-----------------------------------------------------------------------------
\section{The Miller loop}
\label{section:loop}
%-----------------------------------------------------------------------------

The procedure to evaluate the function $f_{n,P}$ is given by Miller in
\cite{Miller04} and relies on the recursive relations
$f_{1, P} = 1$ and
\begin {equation}
\label {eq:millerrecursion}
f_{n+m, P} = f_{n,P} f_{m,P} \frac {\ell_{n P, m P}}{v_{(n + m) P}},
\end {equation}
which can be easily checked by taking divisors on both sides of
the equality, see \S\ref {subsection:divisors}.
This allows to compute $f_{x, P}$ alongside $x P$ in a standard
double-and-add fashion, resulting in the special case of
Algorithm~\ref {algo:millernaf} in which all digits $x_i$ are $0$ or~$1$.
To avoid field inversions, the numerator and the denominator of the
function are kept separate; the correctness of the algorithm may be
derived from the loop invariant
$f_{x',P}(Q) = f/g$, where $x'=\sum_{j=0}^{n-i-1} x_{i+j+1} 2^j$.

\begin{algorithm}[!htb]
\SetKw{KwDownto}{downto}
\SetCommentSty{textit}
\DontPrintSemicolon
\KwIn{%
  $P\ne\mathcal{O}, Q \ne \lambda P$, two points on an elliptic curve $E$ over a field \newline
  $x=\sum_{i=0}^n x_i 2^i$ with $x_i \in \{-1, 0, 1\}$
}
\KwOut{$f_{x,P}(Q)$}
$R \leftarrow P$\;
$f \leftarrow 1$, $g \leftarrow 1$\;
\For{$i \leftarrow n-1$ \KwDownto $0$}{
  $f \leftarrow f^2 \ell_{R,R}(Q)$\;
  $R \leftarrow 2 R$\;
  $g \leftarrow g^2 v_{R}(Q)$\;
  \If{$x_i=1$}{
    $f \leftarrow f \ell_{R,P}(Q)$\;
    $R \leftarrow R+P$\;
    $g \leftarrow g v_{R}(Q)$\;
  }
  \If{$x_i=-1$}{
    $f \leftarrow f \ell_{R,-P}(Q)$\;
    $R \leftarrow R-P$\;
    $g \leftarrow g v_{R}(Q) v_{P}(Q)$
  }
}
\KwRet{$f/g$}
\caption{Miller's algorithm using a NAF}
\label{algo:millernaf}
\end{algorithm}
Techniques for speeding up scalar products on elliptic curves may be adapted.
For instance, allowing addition-subtraction chains for computing $x P$,
preferably by writing $x$ in non-adjacent form (NAF), in which no two
consecutive digits are non-zero, results in Algorithm~\ref {algo:millernaf}.
Notice the additional multiplication by $v_P$ in the case of a digit $-1$,
due to $f_{-1, P} = 1 / v_P$, which is often incorrectly left out in the
literature; in particular, \cite [Algorithm~1]{Miller04} is only
correct for addition chains.
The common omission may be explained by the use of denominator elimination:
If $k$ is even and $G_2 = E (\fqk) \cap \ker (\pi - q)$, then the
$x$-coordinates of $P$ and $Q$ and thus all
values $v_R (Q)$ lie in the subfield $\F_{q^{k/2}}$ and disappear for
reduced pairings involving a final exponentiation. In particular, $g$ may
be dropped completely.

An observation made in \cite{BoElLaLe10} allows to simplify the expression
in the case $x_i = -1$ also for odd~$k$.
Denote by $\lambda$ the slope of the line between $R$ and $-P$. Then
\[
\ell_{R, -P} (Q) = y (Q) + y (P) - \lambda \big( x (Q) - x (P) \big),
\]
which is computed with one multiplication, and $f$ and $g$ are updated
with four multiplications altogether, or with two multiplications if
denominator elimination applies. Then
\begin {equation}
\ell'_{R, -P} (Q) := \frac {\ell_{R, -P} (Q)}{v_P (Q)}
= \frac {y (Q) + y (P)}{x (Q) - x (P)} - \lambda
\label {eq:lprime}
\end {equation}
is obtained without any multiplication since the first term may be
precomputed once and for all. So replacing the block in the case of
$x_i =  -1$ by
\[
f \leftarrow f \ell'_{R, -P} (Q), \quad
R \leftarrow R - P, \quad
g \leftarrow g v_R (Q)
\]
reduces the number of multiplications to~$2$, the same as in the presence
of denominator elimination.

\begin{algorithm}[!htb]
\SetKw{KwDownto}{downto}
\SetCommentSty{textit}
\DontPrintSemicolon
\KwIn{%
  $P\ne\mathcal{O}, Q \ne \lambda P$, two points on an elliptic curve $E$ over a field \newline
  $x=\sum_{i=0}^n x_i 2^i$ with $x_i \in \{-1, 0, 1\}$
}
\KwOut{$f_{x,P}(Q)$}
$h \leftarrow n + \#\{x_i|x_i\ne0\}_{0\le i \le n} - 1$\;
$\delta \leftarrow (-1)^h$\;
$g \leftarrow f_{\delta, P} (Q)$, where $f_{1, P} = 1$ and $f_{-1, P} = 1 / v_P$\;
$f \leftarrow 1, R \leftarrow P$\;
\For{$i \leftarrow n-1$ \KwDownto $0$}{
  \If{$\delta=1$}{
    $f \leftarrow f^2 \ell_{R,R}(Q)$\tcp*[r]{ $f_{-2n, P} = (f_{n,P}^2 \ell_{[n]P, [n]P})^{-1}$}
    $g \leftarrow g^2$\;
    $R \leftarrow 2 R,\ \delta \leftarrow -\delta$\;
    \If(\tcp*[f]{$f_{n+1, P} = (f_{-n,P} f_{-1,P} \ell_{[-n]P, -P})^{-1}$}){$x_i=1$}{
      $g \leftarrow g \ell'_{-R,-P}(Q)$\;
      $R \leftarrow R+P,\ \delta \leftarrow -\delta$\;
    }
    \If(\tcp*[f]{$f_{n-1, P} = (f_{n,P} \ell_{[-n]P, P})^{-1}$}){$x_i=-1$}{
      $g \leftarrow g \ell_{-R,P}(Q)$\;
      $R \leftarrow R-P,\ \delta \leftarrow -\delta$\;
    }
  } \Else {
    $g \leftarrow g^2 \ell_{-R,-R}(Q)$\tcp*[r]{$f_{2n, P} = (f_{-n,P}^2 \ell_{[-n]P, [-n]P})^{-1}$}
    $f \leftarrow f^2$\;
    $R \leftarrow 2 R,\ \delta \leftarrow -\delta$\;
    \If(\tcp*[f]{$f_{-(n+1), P} = (f_{n,P} \ell_{[n]P, P})^{-1}$}){$x_i=1$}{
      $f \leftarrow f \ell_{R,P}(Q)$\;
      $R \leftarrow R+P,\ \delta \leftarrow -\delta$\;
    }
    \If(\tcp*[f]{$f_{-(n-1), P} = (f_{n,P} f_{-1,P} \ell_{[n]P, -P})^{-1}$}){$x_i=-1$}{
      $f \leftarrow f \ell'_{R,-P}(Q)$\;
      $R \leftarrow R-P,\ \delta \leftarrow -\delta$\;
    }
  }
}
\KwRet{$f/g$}
\caption{Boxall \emph{et al.}'s algorithm using a NAF}
\label{algo:caennaf}
\end{algorithm}

In the same article \cite {BoElLaLe10}, Boxall \textit {et al.} introduce
a variant of the algorithm based on
\[
f_{n+m,P}=(f_{-n,P} f_{-m,P} \ell_{-n P, -m P})^{-1},
\]
which contains only three factors instead of the four in
\eqref {eq:millerrecursion}, but requires a swap of the numerator and the
denominator for each doubling or adding on the elliptic curve.
The algorithm of \cite{BoElLaLe10} is given only for addition chains,
but can be generalised to addition-subtraction chains, yielding
Algorithm~\ref {algo:caennaf}.
The correctness of the algorithm stems from the equations given in
commentary and \eqref {eq:lprime}, and the loop invariant
$f_{\delta x',P}(Q)^\delta = f/g$,
where $x'=\sum_{j=0}^{n-i-1} x_{i+j+1} 2^j$.
The value of $h$ is the total number of curve doublings, additions and
subtractions carried out during the algorithm; $\delta$ is $+1$ or $-1$,
respectively, depending on whether the number of curve operations still to
be carried out is odd or even. In the end, $\delta = 1$, and $f/g$ is the
desired result.

%-----------------------------------------------------------------------------
\section{Final exponentiation}
\label{section:exponentiation}
%-----------------------------------------------------------------------------

After the evaluation of a rational function, most pairings require a final
powering to obtain a uniquely defined result. It has been suggested that at
high levels of security the final exponentiation would become so
computationally expensive that the Weil pairing should be preferred to the
Tate pairing \cite{KoMe05}, but this conclusion was quickly
contradicted by a finer analysis of the exponentiation step~\cite{GrPaSm06}.
The rather special form of the exponent $(q^k-1)/r$ makes the final powering
less daunting than it may first appear. This exponent can be decomposed
into three parts as follows. Let $i$ be the smallest prime divisor of $k$.
\begin{equation*}
\frac{q^k-1}{r}
= \underbrace {(q^{k/i}-1)
\cdot \frac{q^k-1}{(q^{k/i}-1) \Phi_k(q)}}_{\text {easy}}
\cdot \underbrace {\frac{\Phi_k(q)}{r}}_{\text {hard}}
\end{equation*}
where $\Phi_k$ is the $k$-th cyclotomic polynomial. The first two factors are
sums of powers of $q$ and are easily computed using only a few applications
of the Frobenius map $x \mapsto x^q$
and multiplications in $G_T$, together with an extra
division in $G_T$ for the first part.
Consequently, the powering by $N = \Phi_k(q)/r$, dubbed the ``hard part'',
is the only truly expensive step of the final exponentiation.
It is usually carried out with mul\-ti-ex\-po\-nen\-ti\-a\-tion
techniques.

\paragraph{Generic exponentiation}
An algorithm proposed by Avanzi and Mih\u{a}ilescu in~\cite{AvMi04}
makes particularly intensive use of the Frobenius map, which is very
efficient for finite field extensions given by binomials, see
Table~\ref {table:opratios}.
To compute $b^N$, it first decomposes $N$ in base~$q$, then each
coefficient in base~$2^\ell$ for a suitably chosen small~$\ell$, such
that
$N = \sum_{i=0}^{\lceil \log_q N \rceil}
\sum_{j = 0}^{\lceil \log_{2^\ell} (q-1) \rceil}
n_{ij} q^i 2^{j \ell}$, and
$b^N = \prod_j \left( \prod_i b^{n_{ij} q^i} \right)^{2^{j \ell}}$.
After precomputing all possible values of $b^{n_{ij}}$ with about
$2^\ell$ multiplications, for each~$j$ the inner product is evaluated
in a Horner scheme; altogether, this requires $O ((\log N) / \ell)$
Frobenius maps and as many multiplications.
The outer product is then computed by a
square-and-multiply approach with $O (\log q) = O ((\log N) / k)$
multiplications, most of which are squarings. Notice that for
$2^\ell \in O (\log N / \log \log N)$ and $k$ growing sufficiently
fast with respect to~$N$ (which is the case due to the security
equivalence of \S\ref {section:curve}), the total complexity
becomes $O (\log N / \log \log N)$ operations, which is sublinear
in the bit size of~$N$.

As can be seen in Table~\ref{table:opratios}, the Frobenius is not
completely negligible when the field extension is realised by a trinomial,
so we investigated an alternative approach due to Nogami \emph {et al.}
\cite{NoKaNeMo08}, which purportedly requires fewer applications of the
Frobenius map.
Let $\ell$ be an integer and $c = \lceil (\log_q N) / \ell \rceil$,
and let $t = \lceil \log_2 q \rceil$.
The algorithm of \cite{NoKaNeMo08} creates a binary matrix with
$\ell$ rows and $c t$ columns by first writing $N$ in basis $q^c$.
Each coefficient corresponds to one row and is decomposed into $c$
coefficients in base $q$, each of which is in turn written with $t$
coefficients in base $2$. These form the columns of the matrix,
organised into $c$ blocks of $t$ columns each.
To compute $b^N$, first the powers $b^{2^i}$ are precomputed with
$t-1$ multiplications. If the same column occurs $d \geq 2$ times
in $e \leq c$ column blocks, its occurrences can be combined with $d-1$
multiplications and $e-1$ applications of the Frobenius map.
Taking into account that there are at most $2^\ell - 1$ different
non-zero columns, this step can thus be carried out with at most
$c t - 1$ multiplications and $(c-1) (2^\ell - 1)$ Frobenius maps.
Heuristically, a fraction of $1/2^\ell$ of the columns are zero, for which
there is nothing to do; so one expects a number of multiplications
closer to $c t (1 - 1/2^\ell)$, a noticeable difference for small values of
$\ell$.
For each row, the combined columns which have a $1$ in this row need to
be multiplied, with at most $\ell (2^{\ell-1} - 1)$ multiplications.
Finally, the rows are combined with $\ell-1$ multiplications and $\ell-1$
Frobenius maps. We arrive at a complexity of
$(c - 1) (2^\ell - 1) + \ell - 1$ Frobenius maps and
$c t + \ell (2^{\ell-1} - 1) + \ell + t - 3$ multiplications.
The latter can be tightened heuristically by multiplying the first
term with $1 - 1/2^\ell$, which turns out to be close to the experimentally
observed values.
Asymptotically, the number of multiplications is equivalent to
$(\log_2 N) / \ell + \ell 2^\ell$.
Recall from the security equivalence of~\S\ref {section:curve}
that $k$ is of the order of $\log^2 q$, so that
$c \in O \left( (\log N)^{2/3} / \ell \right)$
and the number of Frobenius maps is in
$O \left( 2^\ell (\log N)^{2/3} / \ell \right)$.
Letting $2^\ell = (\log N)^{1/3}$ yields a sublinear total complexity
of $O (\log N / \log \log N)$.
The analysis also shows that by preferring smaller values of~$\ell$, one may
reduce the number of Frobenius maps compared to Avanzi and Mih\u{a}ilescu's
algorithm, at the price of more multiplications.
Notice that $\ell$ fixes~$c$, and that the exponentiality in~$\ell$ implies
that an optimal value is found easily in practice.

\paragraph{Family-dependant exponentiation}
Scott \emph{et al.} proposed in \cite{ScBeChPeKa09} an exponentiation technique
for polynomial families of parameters $q (x)$ and $r (x)$.
The exponent is written first in base $q (x)$, then each coefficient
in base $x$ as \\
\centerline {
$N (x) = \Phi_k (q (x)) / r (x)
= \sum_{i = 0}^{\lfloor \deg N (x) / \deg q (x) \rfloor}
  \sum_{j = 0}^{\deg q(x) - 1}
  \lambda_{i, j} x^j q (x)^i$.}
To obtain $b^{N (x_0)}$, the values $b^{x_0^j q^i}$ are precomputed with
first $\deg q (x)$ exponentiations by $x_0$,
which can be done with $O (\log q)$ multiplications, then
about $\deg N (x)$ applications of the Frobenius map.
The final result is then obtained by multi-exponentiation with the exponents
$\lambda_{i, j}$. The exact complexity of this step depends on the length
of an addition sequence passing through all of the $\lambda_{i, j}$.
If $\Lambda = \max |\lambda_{i, j}|$, then the best rigorous bound currently
available is
$\log_2 \Lambda + \deg N (x) O (\log \Lambda / \log \log \Lambda)$,
where $\deg N (x)$ is the number of potentially different values of
$\lambda_{i, j}$, realised, for instance, by~\cite {Yao76}.
In practice,
the coefficients are small, and there are addition sequences with only few
additional terms, see Table~\ref {table:addseq},
leading to a heuristic complexity of $O (\deg N (x))$ multiplications.
The total complexity then becomes $O (\log q + \deg N (x))$ multiplications
and $O (\deg N (x))$ Frobenius maps,
where $\deg N (x) \approx \log N / \log x_0$.

%-----------------------------------------------------------------------------
\section{Implementation}
\label{section:implementation}

We have implemented the various pairings for the different curves of
\S\ref {section:curve} in the PARI library and linked them into the free
number theoretic computer algebra system GP \cite{parigp}. Our aim was not
to provide an optimal \textit {ad hoc} implementation for any one of the
curves or pairings, but rather to keep a sufficient level of genericity
appropriate for a general purpose system, while still implementing
algorithmic optimisations that apply in a broader context.
All benchmarks were performed on a Macbook Pro with a 2.5 GHz Core 2 Duo
processor, and timings are given in milliseconds (rounded to two
significant digits).

\subsection{Finite field arithmetic}
\label{ssec:ff}

Arithmetic in $\fp$ is already available in PARI. The field extensions
$\fpk$ are defined by binomials $X^k-a$ for all curves but \curveg, \curveh\
and \curveo, for which only trinomials of the form $X^k + X + a$
exist. A defining binomial can be found if and only if
all prime factors of~$k$ divide $p-1$, and, additionally for $4 \mid k$,
if $p \equiv 1 \pmod 4$ \cite [Theorem~2]{PaTh09}, which happens quite often
in our context where $k$ has only few prime factors and there is a certain
freedom in the choice of~$p$.
Fast reduction modulo binomials and trinomials had to be added to PARI.
Definition by a binomial is crucial for an efficient pairing implementation,
see Table~\ref{table:notwist}, which shows the effect for the ate pairing
on \curvei\ when artificially switching to a trinomial.
\begin{table}[!hbt]
\centerline {
\begin{tabular}{c|c|c}
	Defining polynomial & unreduced ate & unreduced optimal \\
	\hline
  $x^{18} + x + 6$    & 490 ms        & 120 ms\\
	$x^{18} + 19$       & 150 ms        &  35 ms
	\end{tabular}
}
  \caption {Timings for ate pairing depending on the finite field defining polynomial}
  \label{table:notwist}
\end {table}

In several places in the literature, it is suggested to build $\fpk$
by successive extensions of degree dividing~$k$, in particular of
degree~$2$ or~$3$. It is unclear where the benefits of this strategy lie
for multiplication: Virtually the same effect may be reached by using
Karatsuba (degree~$2$) and Toom-Cook (degree~$3$) for polynomial arithmetic,
which moreover speeds up the computations also when $k$ is not a power
product of~$2$ and~$3$. By keeping a single extension, it is also easier
to reach the thresholds for asymptotically fast algorithms.
In particular, PARI uses Kronecker substitution to replace the product of
polynomials by that of large integers which is outsourced to the GMP
library \cite{gmp}; in experiments, this turned out to be faster than
Karatsuba multiplication.

Note that using a binomial to define the field extensions also indirectly
speeds up the arithmetic when field elements lie in subfields of $\fqk$,
which happens systematically for the group $G_2 = E [r] \cap \Ker (\pi - q)$
in the presence of twists.
As an example consider again the curve \curvei.
Let $D \in \F_{p^3}$ be a quadratic and a cubic non-residue in $\F_{p^3}$,
which implies that $\F_{p^{18}} = \F_{p^3}[D^{1/6}]$.
Then \curvei$: y^2 = x^3 + b$ admits a sextic twist
$E' : y^2 = x^3 + b/D$, explicitly given by
$\phi_6 : E' (\F_3) \rightarrow E(\F_{18}) :
(x,y) \mapsto \left( \sqrt[3]{D} x, \sqrt{D} y \right)$,
which yields an isomorphism of $E' (\F_{p^3})[r]$ with~$G_2$.
If $\F_{18}$ is defined by a binomial $X^{18} + A$, then $D = A^{1/3} = X^6$,
the elements of $\F_{p^3}$ are written as $a_{12} X^{12} + a_6 X^6 + a_0$,
and an element $Q$ of $G_2$ is given as
$\phi_6 (Q') = (a_{14} X^{14} + a_8 X^8 + a_2 X^2, a_{15} X^{15} +
a_9 X^9 + a_3 X^3)$ with all $a_i \in \F_p$.
These sparse polynomials lead to a faster arithmetic,
and part of the speed gains for binomial field extensions as shown in
Table~\ref {table:notwist} may be attributed to this implicit handling
of subfields.

Explicit towers of finite fields could be useful for realising the Frobenius
automorphism of $\fpk/\fp$, since the non-trivial automorphism of a
quadratic extension is a simple sign flip; while those of
cubic extensions require a multiplication in the smaller field anyway.
We chose instead to implement the Frobenius, as well as its powers, as
linear maps by multiplying with precomputed $k\times k$ matrices.
We originally intended to study the use of optimal normal bases of
$\fpk / \fp$, in which the Frobenius $\pi$ is realised as a simple permutation
of the coordinates. It turns out, however, that again binomial field extensions
yield an efficient arithmetic: In their presence, the matrix of $\pi$
is the product of a diagonal matrix and a permutation,
so $\pi$ and its powers can be computed with~$k$ multiplications
in~$\fp$ \cite[Theorem~3]{BaPa01}.
In the trinomial case, the Frobenius is computed
with close to $k^2$ multiplications in $\fp$.

Table~\ref {table:opratios} summarises the relative costs of the Frobenius~$\pi$,
the multiplication $M_2$ and the inversion $I_2$ in $\fpk$ and the
multiplication $M_1$ and the inversion $I_1$ in $\fp$.
The effect of a defining trinomial on the cost of $\pi$ is clearly visible.
\begin{table}[!hbt]
\centerline{
\scriptsize
  \begin{tabular}{c|r|r|r||r|r|r|r|r|r||r|r|r|r|r}
               & \curvez & \curveg & \curvefii & \curvefiv & \curvep & \curves & \curvet & \curvei & \curveh & \curvefv &  \curvel & \curveo & \curven & \curvem \\
\hline
    $M_1 / \mu$s    &   0.41  &  0.36   &     0.27  &    0.82   &   0.64  &  0.57   &  0.66   &   0.64  &   0.49  &   2.2    &      0.74 &  0.83  &   0.73  &   0.64  \\
    $I_1/M_1$       &   11    &  11     &  15       &  12       &  13     &  14     &  12     &  12     &  13     &  10      &   13     &  12     &  12     &  13     \\
    $M_2/M_1$       &   55    &  90     & 110       &  80       & 120     &  130    &  130    & 140     & 170     &  70      &  250     & 210     & 250     & 240     \\
    $I_2/M_2$       &    8.0  &   7.9   &   8.6     &   8.1     &   8.1   &  8.7    &  9.2    &   8.8   &   8.7   &   8.1    &    9.2   &   9.5   &  10     &  10     \\
    $\pi/M_2$       &   0.19  &   0.63  &   0.19    &   0.17    &   0.15  & 0.15    & 0.15    &   0.16  &   0.95  &   0.18   &    0.14  &   1.2   &   0.16  &   0.15
\end{tabular}
}
\caption{Costs of arithmetic in finite prime and extension fields}
  \label{table:opratios}
\end {table}

\subsection{Miller loop}
\label {ssec:miller}

Given the cost of inversion in Table~\ref {table:opratios}, we implemented
the elliptic curve arithmetic using affine coordinates. Timings for the Miller
loop are summarised in Table~\ref {table:millertimes}. While mathematically
not necessary, the Tate and Weil pairings have also been restricted to the
subgroups $G_1$ and $G_2$ of eigenvalue~$1$ and~$p$, respectively, which
assures a fairer comparison and incidentally a type~$3$ pairing in the
notation of \cite{GaPaSm08}, see also \cite{ChMe11}.
For even embedding
degree, we applied denominator elimination. The first row uses a
double-and-add approach, the second one a signed NAF with
Algorithm~\ref {algo:millernaf}. The variant of
Algorithm~\ref {algo:caennaf} is only of interest when $k$ is odd; we give
its timings with a double-and-add chain and a NAF in the third and fourth
rows, respectively.
It makes an impressive difference.
\begin{table}[!hbt]
\centerline {
\begin{tabular}{c|r|r|r|r|r|r|r|r|r|r}
	  Curve      & $\et$ & $\etw$  & $\evt$ &   $\ea$  & $\ev$  & $\ew$ & $\deg q$ & $\deg r$ & $\deg t$ & $\varphi(k)$\\
	  \hline
    \curvez    &  31   &  64     &    15  &       50 &    14  &   100 &        8 &        6 &        4 &   6 \\
               &  29   &  59     &    15  &       45 &    14  &    97 &          &          &          &     \\
               &  21   &  46     &    11  &       46 &    12  &    85 &          &          &          &     \\
               &  20   &  42     &    11  &       42 &    12  &    80 &          &          &          &     \\
    \hline                                                                                        
    \curveg    &  43   &    --   &     -- &      110 &     20 &  230  &       24 &       20 &       12 &  10 \\
               &  39   &    --   &     -- &      100 &     20 &  200  &       24 &       20 &       12 &  10 \\
               &  31   &    --   &     -- &      100 &     20 &  200  &          &          &          &     \\
               &  28   &    --   &     -- &       96 &     19 &  180  &          &          &          &     \\
    \hline                                                                                 
    \curvefii  &  14   &    13   &      7 &       18 &      9 &    70 &        4 &        4 &        2 &   4 \\
               &  14   &    13   &      7 &       18 &      9 &    67 &          &          &          &     \\
    \hline                                                                                 
    \hline                                                                                 
    \curvefiv  & 93    &    91   &     53 &      110 &     54 &   420 &        4 &        4 &        2 &   4 \\
               & 91    &    91   &     52 &      110 &     54 &   410 &          &          &          &     \\
    \hline                                                                                 
    \curvep    & 130   &   520   &     73 &      270 &     44 &   460 &       12 &        8 &        6 &   8 \\
               & 130   &   480   &     71 &      240 &     41 &   440 &          &          &          &     \\
               &  90   &   360   &     53 &      240 &     41 &   370 &          &          &          &     \\
               &  88   &   330   &     50 &      220 &     37 &   360 &          &          &          &     \\
    \hline                                                                                 
    \curves    &  64   &  170    &     35 &      150 &    230 &  360  &       10 &        8 &        5 &   8 \\
               &  62   &  160    &     34 &      140 &    220 &  340  &          &          &          &     \\
    \hline                                                                                 
    \curvet    &  160  &   --    &     -- &       54 &    54  & 930   &       38 &        32 &        2 &   16 \\
               &  140  &   --    &     -- &       50 &    50  & 820   &          &          &          &     \\
               &  110  &   --    &     -- &       51 &    51  & 830   &          &          &          &     \\
               &  100  &   --    &     -- &       47 &    47  & 740   &          &          &          &     \\
    \hline                                                                                 
    \curvei    &  78   &   160   &     39 &      160 &     35 &  400  &        8 &        6 &        4 &   6 \\
               &  75   &   150   &     36 &      150 &     35 &  380  &          &          &          &     \\
    \hline                                                                                 
    \curveh    & 130   &    --   &     -- &      350 &     40 &   760 &       40 &       36 &       20 &  18 \\
               & 120   &    --   &     -- &      320 &     35 &   680 &          &          &          &     \\
               &  93   &    --   &     -- &      340 &     40 &   670 &          &          &          &     \\
               &  88   &    --   &     -- &      300 &     34 &   610 &          &          &          &     \\
    \hline                                                                                 
    \hline                                                                                 
    \curvefv   & 410   &   400   &    240 &      450 &    220 &  1800 &        4 &        4 &        2 &   4 \\
               & 400   &   400   &    240 &      450 &    230 &  1700 &          &          &          &     \\
    \hline                                                                                 
    \curvel    & 190   &    88   &     88 &       55 &     55 &   960 &       10 &        8 &        1 &   8 \\
               & 180   &    85   &     85 &       56 &     56 &   900 &          &          &          &     \\
    \hline                                                                                 
    \curveo    & 450   &    --   &     -- &     1400 &    120 &  2700 &       52 &       40 &       26 &  20 \\
               & 410   &    --   &     -- &     1300 &    130 &  2400 &          &          &          &     \\
               & 310   &    --   &     -- &     1400 &    130 &  2400 &          &          &          &     \\
               & 300   &    --   &     -- &     1200 &    120 &  2200 &          &          &          &     \\
    \hline                                                       
    \curven    & 210   &    --   &     -- &      660 &     99 &  1600 &       28 &       24 &       14 &  12 \\
               & 190   &    --   &     -- &      610 &     90 &  1400 &          &          &          &     \\
    \hline                                                       
    \curvem    & 370   &  1900   &    190 &      570 &     54 &  1300 &       20 &       18 &       10 &  18 \\
               & 350   &  1800   &    180 &      510 &     53 &  1300 &          &          &          &     \\
               & 260   &  1300   &    130 &      530 &     52 &  1100 &          &          &          &     \\
               & 250   &  1200   &    130 &      450 &     47 &  1100 &          &          &          & 
	\end{tabular}
}
	\caption{Timings of the function evaluation in milliseconds}
  \label{table:millertimes}
\end {table}

Generically, one expects a NAF to save about 11\% of the number of
operations. For our curves, the effect is often much less. This can be
explained by the sparsity of the integer~$r$ derived from a curve family,
which is thus closer to non-adjacent form than a random integer.
For instance, the binary decomposition of~$r$ for \curvefii\ has
only $87$ entries $1$ out of $256$, a density that would be expected
in a NAF of a random number. The NAF has $37$ entries $1$ and $-1$ each.
Also counting the squarings, the gain in the number of operations is
less than~4\%.
One could reduce the number of multiplications even further by combining
with a sliding window technique; since the number of squarings is
unchanged, the effect will be more and more marginal with an increasing
window size.

At these high security levels and consequently big values of~$k$, the ate
pairing~$\ea$ is clearly in general slower than the Tate pairing; the
smaller multiplier is more than offset by the need to work over the extension
field. The only exception is \curvel\ with a particularly small trace
($t(x)$ of degree~$1$ for $r(x)$ of degree~$8$). In fact, the ate pairing for
this curve is already optimal.
As can be expected, the twisted ate pairing $\etw$ is even less efficient
except for small values of~$k$ combined with a high degree twist:
The power of the trace $T^e = T^{k/d}$ quickly exceeds $r$ itself
(in the table, we computed with $T^e$; one could reduce modulo~$r$ and
arrive at the same timings as the Tate pairing, but may then as well
stick with the original).
The optimal versions indeed keep their promises. Due to the overhead of
computing several functions, the total running time is not reduced by
a factor of~$\varphi (k)$, but the optimal ate pairing is generally
faster than the Tate pairing.
Twisted pairings are asymptotically slower, but interesting for the
medium values of~$k = 9$, $12$ or $16$ which admit a twist of (relatively)
high degree $3$, $6$ or~$4$.

The Weil pairing with its two function evaluations could be expected to
behave like a Tate followed by an ate pairing; due to the different loop
lengths, the part $f_{r, Q} (P)$, however, has a complexity closer to
$\deg r (x) / \deg t (x)$ times that of ate, as can be roughly verified in
the table.
As already stated in the literature, the enormous overhead of the Weil
pairing is not offset by saving the final exponentiation, see
\S\ref {subsection:exponentiation}.

At higher security levels, odd values of~$k$ lead to a bigger $\varphi (k)$
and thus a higher gain in the optimal pairing; together with the
Miller loop improvement of~\cite {BoElLaLe10}, odd and in particular prime
values of~$k$ become attractive. Notice that \curveg, \curveh\ and
\curveo\ are heavily penalised by the trinomial field extension. Indeed,
odd or, worse, prime values of~$k$ make the divisibility conditions for the
existence of a binomial extension harder to satisfy. Moreover, the degree of
$p(x)$ also grows with $\varphi (k)$, so that the polynomial represents
fewer numbers in the desired range and leaves less choice for~$p$ or
a value of $x_0$ with low Hamming weight.
Even if a binomial field extension exists, odd values of~$k$ that are not
divisible by~$3$ (in particular, prime~$k$ again) suffer from a lack of
twists and thus a less efficient field arithmetic as discussed in
\S\ref {ssec:ff}.

\subsection{Final exponentiation}
\label {subsection:exponentiation}

Timings for the final exponentiation are compiled in
Table~\ref {table:exptimes}.
\begin{table}[!hbt]
\centerline {
\begin{tabular}{c|r|r|r|r|r|r}
	  Curve       &   $\varphi (k)$ & Naive & Hard naive & AM04 & NMKM08 & SBCPK09\\
	  \hline                                   
    \curvez     &   6 &      56  &      36  &     15 &      15 &     8 \\
    \hline                    
    \curveg     &  10 &      80 &       78   &    28  &     24 &    21 \\
    \hline
    \curvefii   &   4 &      58  &      17   &     8  &      8 &     4 \\
    \hline
    \hline
    \curvefiv   &   4 &     380  &    100   &    41  &     37 &    26 \\
    \hline
    \curvep     &   8 &     490  &    250   &    85  &    110 &    50 \\
    \hline
    \curves     &   8 &     420  &    200   &    66  &     80 &    44 \\
    \hline
    \curvet     &  16 &     580  &    550   &   180  &    200 &   110 \\
    \hline
    \curvei     &   6 &     680  &    210   &    78  &     83 &    49 \\
    \hline
    \curveh     &  18 &     460  &    460   &   150  &    110 &    83 \\
    \hline
    \hline
    \curvefv    &   4 &    1800  &    470   &   170  &    170 &   120 \\
    \hline
    \curvel     &   8 &    2000  &    640   &   150  &    200 &    97 \\
    \hline
    \curveo     &  20 &    2600  &   2100   &   700  &    470 &   320 \\
    \hline
    \curven     &  12 &    2300  &   1000   &   240  &    270 &   170 \\	
    \hline
    \curvem     &  18 &    2100  &   1400   &   290  &    310 &   130
\end{tabular}
}
	\caption{Final exponentiation times in milliseconds}
  \label{table:exptimes}
\end {table}

The first column corresponds to a direct exponentiation by $(q^k-1)/r$
via the sliding window algorithm built into PARI. The second column
does so for the hard part, while computing the easy one using
Frobenius maps. The next two columns relate the implementation of the
hard parts following~\cite{AvMi04} and~\cite{NoKaNeMo08}.
At low security level, the differences between these two algorithms
are minimal. For the medium and high security range, our implementation
confirms the claim of~\cite{NoKaNeMo08}: Their algorithm becomes faster
when Frobenius maps are more expensive, as for the three curves
\curveg, \curveh\ and \curveo.
The $k=12$ curves stand out: The low value of $\varphi (k)$ makes the
final exponentiation much easier with these two algorithms that rely
on an expansion in base $q$.

While the theoretical analysis of \S\ref {section:exponentiation}
is not conclusive, the experiments are unequivocal:
The algorithm of \cite{ScBeChPeKa09} (which we used, as explained in the
article, to potentially compute a small power of the true pairing if the
coefficients of the polynomial contain denominators) is clearly the
fastest one for curves obtained from polynomial families.

%-----------------------------------------------------------------------------
\section{Overall timings and conclusion}
\label{section:conclusion}
%-----------------------------------------------------------------------------

For each of our reference curves, Table \ref{table:fastest} summarises the
timings obtained for the fastest pairing.

\begin{table}[!hbt]
\centerline {
	\begin{tabular}{l|c|r|c|r|r|r}
Security & Curve     & $\varphi (k)$ & Pairing   &  Unreduced   & Final exp & Reduced \\
	  \hline
128 bit & \curvez    &           6   & $\evt$    & 11           & 8         &  19 \\
\cline {2-7}                                   
        & \curveg    &           10  & $\ev$     & 19           & 21        &  40 \\
\cline {2-7}                                   
        & \curvefii  &           4   & $\evt$    &  7           & 4         &  11 \\
    \hline                                     
    \hline                                     
192 bit & \curvefiv  &           4   & $\evt$    & 52           & 26        &  78 \\
\cline {2-7}                                     
        & \curvep    &           8   & $\ev$     & 37           & 50        &  87 \\
\cline {2-7}                                     
        & \curves    &           8   & $\evt$     & 34           & 44        & 78 \\
\cline {2-7}                                     
        & \curvet    &          16   & $\ea=\ev$ & 47           & 110       & 157 \\
\cline {2-7}
        & \curvei    &           6   & $\ev$     & 35           & 49        &  84 \\
\cline {2-7}                                     
        & \curveh    &          18   & $\ev$     & 34           & 83        & 120 \\
    \hline                                       
    \hline                                       
256 bit & \curvefv   &           4   & $\ev$     & 220          & 120       & 340 \\
\cline {2-7}                                     
        & \curvel    &           8   & $\ea=\ev$ & 55           & 97        & 150 \\
\cline {2-7}                                     
        & \curveo    &          20   & $\ev$     & 120          & 320       & 440 \\
\cline {2-7}                                     
        & \curven    &          12   & $\ev$     & 90           & 170       & 260 \\
\cline {2-7}                                     
        & \curvem    &          18   & $\ev$     & 47           & 130       & 180
	\end{tabular}
}
	\caption{Timings of the fastest reduced pairing variants}
  \label{table:fastest}
\end {table}

Optimal pairings are indeed optimal for higher security levels.
Their unreduced version benefits from high values of $\varphi (k)$,
as can be seen by comparing \curvel\ and \curvem. However, part of this
advantage is offset by the lack of denominator elimination for odd~$k$,
although Boxall \textit {et al.}'s variant almost closes the gap again.
Moreover, the higher cost for the final exponentiation more than compensates
the gain in the Miller loop.
The decision which pairing to take then also depends on the concrete
cryptographic protocol: Not all of them require reduced pairings throughout
their execution. For instance, verification protocols such
as~\cite {BlFuIzJaSeVe10} make do with testing equality of products of
several pairings. All of these may then be computed unreduced, and only a
final quotient of products needs to be raised to the power, which makes
this exponentiation negliglible.

For a reduced pairing at lower security levels, Barreto--Naehrig
curves with $k=12$ remain unbeaten, profiting from an exceptionally fast
final exponentiation.

At 192 bit security, Barreto--Naehrig curves need to work with a larger
than optimal size of the underlying elliptic curve, but still provide the
fastest pairings.
An equivalent performance, however, may be reached for $k=16$ with
curve~\curves. The suboptimal $\rho = 5/4$ notwithstanding, this curve is
of size 501~bits instead of 663~bits for the Barreto--Naehrig curve,
resulting in less bandwidth for exchanging curve points.
Thus our study shows that \curves\ is preferable at medium security level.

At the highest AES equivalent of 256 bit, Barreto--Naehrig curves are no
longer competitive speed-wise. Here the curve \curvel\ stands out. Although
gaining only a factor of~$8$ in the Miller loop length, it profits from
a very fast final exponentiation, while even the unreduced variant remains
comparable to the closest competitor \curvem.

As becomes clear from this study, extension fields $\fpk$ that do not allow
a binomial as a defining polynomial are to be banned, see \curveg, \curveh\
and \curveo.

Whether odd or even embedding degrees are preferable remains
undecided. Our results seem to indicate that odd degrees are slightly
slower. This can be explained by their higher probability of
requiring a trinomial field extension, sparser families and the lack
of twists as explained
at the end of \S\ref {ssec:miller}. Often, the gain odd and, in particular,
prime embedding degrees provide through larger values of $\varphi (k)$
for the optimal pairings is more than offset by an expensive final
exponentiation, as is well illustrated by~\curvet.
In protocols that work with mostly unreduced pairings, however,
Boxall \textit {et al.}'s variant of the Miller loop makes odd embedding
degrees competitive,
see \curvep\ and \curvem.

A definite conclusion is made difficult by the lack of choice for any given
security level: Some families are so sparse that they contain no curves of
prime cardinality in the desired range or, if they do, no curves allowing
to work with extension fields defined by binomials.
Even if suitable curves exist, the sparsity of a family may have a big
impact on the efficiency of the Miller loop. Notice that the loop
for \curvem\ is more than twice shorter than that of \curvel. Nevertheless,
the unreduced pairing is computed in almost the same time. This can be
explained by the Hamming weight of the multiplier: The family of \curvem\
is instantiated with $x_0$ of weight~$13$, that of \curvel\ with $x_0$ of
weight~$7$.
So the search for new curve families remains a research topic of interest,
not only for families with optimal~$\rho$, as witnessed by
the good performance of \curvep\ despite its very bad $\rho = 3/2$.

\appendix
%---------------------------------------------------------------------------
\section{Curve parameters}
\label{appendix:params}
%---------------------------------------------------------------------------
Tables~\ref {table:curves128} to~\ref {table:curves256} give the exact
parameters for the curves we studied for different security levels.
The following notations are used:
$p(x)$, $r(x)$ and $t (x)$ are the polynomials representing the cardinality
of the finite prime field~$\fq$, a (large) prime factor of the curve
cardinality and the trace of Frobenius, respectively; $x_0$ is the numeric
value of the variable~$x$; $(a,b)$ gives the equation of the curve
$y^2=x^3+ax+b$; $\dpol$ is the irreducible polynomial defining $\fpk$.
For $k=16$ or $17$, no prime values $r (x_0)$ exist in the desired range.
We thus admit a small cofactor and let $r_0$
denote the actual large prime factor of $r (x_0)$.

Table~\ref {table:optivec} provides the short lattice vectors yielding our
optimal ate pairings $\ev$ and optimal twisted ate pairings $\evt$,
see \S\ref {ssec:optimal}.

Table~\ref {table:addseq} records the addition sequences used in the
final exponentiation of \cite{ScBeChPeKa09}. To remove denominators,
the power~$s$ of the original pairing is computed; $n$ is the number of
(not necessarily distinct) non-zero coefficients $\lambda_{ij}$,
see \S\ref {section:exponentiation}.
The underlined terms are those that are added to the sequence.
As can be seen, there is in general a very small number of very small
distinct coefficients, and only a tiny number of terms, if any,
needs to be added.

\begin {table}[!hbt]
\centerline {
\scriptsize
\begin{tabular}{l|l@{\quad}l@{\quad}l}
\curvez
& \multicolumn{3}{l}{$p(x) = (x^8-x^7+x^6-x^5-2x^4-x^3+x^2+2x+1)/3$}\\
& \multicolumn{3}{l}{$r(x) = (x^6-x^3+1)/3$}\\
& \multicolumn{3}{l}{$t(x) = -x^4+x+1$}\\
& {$(a, b) = (0,7)$} & {$\dpol(X)=X^{9}+3$}
& {$x_0 = 43980465324080$}\\
\hline
\curveg
& \multicolumn{3}{l}{$p(x) = (x^{24}-x^{23}+x^{22}-x^{13}+4x^{12}-x^{11}+x^2-x+1)/3$}\\
& \multicolumn{3}{l}{$r(x) = x^{20}+x^{19}-x^{17}-x^{16}+x^{14}+x^{13}-x^{11}-x^{10}-x^9+x^7+x^6-x^4-x^3+x+1$}\\
& \multicolumn{3}{l}{$t(x) = x^{12}+1$}\\
& {$(a, b) = (0,4)$} & {$\dpol(X)=X^{11}+X+11$}
& {$x_0 = 11210$}\\
\hline
\curvefii
& \multicolumn{3}{l}{$p(x) = 36x^4+36x^3+24x^2+6x+1$}\\
& \multicolumn{3}{l}{$r(x) = 36x^4+36x^3+18x^2+6x+1$}\\
& \multicolumn{3}{l}{$t(x) = 6x^2+1$}\\
& {$(a, b) = (0,5)$} & {$\dpol(X)=X^{12}+5$}
& {$x_0 = 6917529027641094616$}
\end{tabular}
}
\caption {\label {table:curves128}
Curves for security level 128 bit}
\end {table}
%------------------------------------------------------------------------------
\begin {table}
\centerline {
\scriptsize
\begin{tabular}{l|l@{\quad}l@{\quad}l}
\curvep
& \multicolumn{3}{l}{$p(x) = (x^{12} - x^{11} + x^{10} - x^7 - 2x^6 - x^5 + x^2 + 2x + 1)/3$}\\
& \multicolumn{3}{l}{$r(x) = x^8 + x^7 - x^5 - x^4 - x^3 + x + 1$}\\
& \multicolumn{3}{l}{$t(x) = -x^6 + x + 1$}\\
& {$(a, b) = (0,13)$} & {$\dpol(X)=X^{15}+13$}
& {$x_0 = 271533021386417$}\\
\hline
\curves
& \multicolumn{3}{l}{$p(x) = (x^{10} + 2x^9 + 5x^8 + 48x^6 + 152x^5 + 240x^4 + 625x^2 + 2398x + 3125)/980$}\\
& \multicolumn{3}{l}{$r(x) = (x^8 + 48x^4 + 625)$\hspace{20.5mm}$r_0 = r(x_0)/20641250$}\\
& \multicolumn{3}{l}{$t(x) = (2x^5 + 41x + 35)/35$}\\
& {$(a, b) = (1,0)$} & {$\dpol(X)=X^{16}+2$}
& {$x_0 = 2251799888961585$}\\
\hline
\curvet
& \multicolumn{3}{l}{$p(x) = (x^{38} + 2x^{36} + x^{34} + x^4 - 2x^2  + 1)/4$}\\
& \multicolumn{3}{l}{$r(x) = x^{32} - x^{30} + x^{28} - x^{26} + x^{24} - x^{22} + x^{20} - x^{18} + x^{16} -x^{14} + x^{12} - x^{10}$}\\
& \multicolumn{3}{l}{\hspace{0.98cm}$ + x^8 - x^6 + x^4 - x^2 + 1$\hspace{12.5mm}$r_0=r(x_0)/12071636373225929$}\\
& \multicolumn{3}{l}{$t(x) = -x^2 + 1$}\\
& {$(a, b) = (13, 0)$} & {$\dpol(X)=X^{17}+2$}
& {$x_0 = 12681$}\\
\hline
\curvei
& \multicolumn{3}{l}{$p(x) = (x^8 + 5x^7 + 7x^6 + 37x^5 + 188x^4 + 259x^3 + 343x^2 + 1763x + 2401)/21$}\\
& \multicolumn{3}{l}{$r(x) = (x^6 + 37x^3 + 343)/343$}\\
& \multicolumn{3}{l}{$t(x) = (x^4 + 16x + 7)/7$}\\
& {$(a, b) = (0,19)$} & {$\dpol(X)=X^{18}+19$}
& {$x_0 = 48422703193491756920$}\\
\hline
\curveh
& \multicolumn{3}{l}{$p(x) = (x^{40}-x^{39}+x^{38}-x^{21}-2x^{20}-x^{19}+x^2+2x+1)/3$}\\
& \multicolumn{3}{l}{$r(x) = x^{36}+x^{35}-x^{33}-x^{32}+x^{30}+x^{29}-x^{27}-x^{26}+x^{24}+x^{23}-x^{21}-x^{20}+x^{18}$}\\
& \multicolumn{3}{l}{\hspace{0.98cm}$-x^{16}-x^{15}+x^{13}+x^{12}-x^{10}-x^9+x^7+x^6-x^4-x^3+x+1$}\\
& \multicolumn{3}{l}{$t(x) = -x^{20}+x+1$}\\
& {$(a, b) = (0,9)$} & {$\dpol(X)=X^{19}+X+23$}
& {$x_0 = 1274$}\\
\hline
\curvefiv
& \multicolumn{3}{l}{$p(x) = 36x^4+36x^3+24x^2+6x+1$}\\
& \multicolumn{3}{l}{$r(x) = 36x^4+36x^3+18x^2+6x+1$}\\
& \multicolumn{3}{l}{$t(x) = 6x^2+1$}\\
& {$(a, b) = (0,13)$} & {$\dpol(X)=X^{12}+5$} & \\
& \multicolumn{3}{l}{$x_0 = 29230032746618058364073696654325660393118650866996$}
\end {tabular}
}
\caption {\label {table:curves192}
Curves for security level 192 bit}
\end {table}

%------------------------------------------------------------------------------
\begin {table}
\centerline {
\scriptsize
\begin{tabular}{l|l@{\quad}l@{\quad}l}
\curvel
& \multicolumn{3}{l}{$p(x) = (x^{10} - 2x^9 + x^8 - x^6 + 2x^5 - x^4 + x^2 + x + 1)/3$}\\
& \multicolumn{3}{l}{$r(x) = x^8-x^4+1$}\\
& \multicolumn{3}{l}{$t(x) = x+1$}\\
& {$(a, b) = (0,1)$} & {$\dpol(X)=X^{24}+19$}
& {$x_0 = 18446744073709602433$}\\
\hline
\curveo
& \multicolumn{3}{l}{$p(x) = (x^{52} - x^{51} + x^{50} - x^{27} - 2x^{26} - x^{25} + x^2 + 2x + 1)/3$}\\
& \multicolumn{3}{l}{$r(x) = x^{40} + x^{35} - x^{25} - x^{20} - x^{15} + x^5 + 1$}\\
& \multicolumn{3}{l}{$t(x) = -x^{26} + x + 1$}\\
& {$(a, b) = (0,31)$} & {$\dpol(X)=X^{25}+X+19$}
& {$x_0 = 6995$}\\
\hline
\curven
& \multicolumn{3}{l}{$p(x) = (x^{28} + x^{27} + x^{26} - x^{15} + 2x^{14} - x^{13} + x^2 - 2x + 1)/3$}\\
& \multicolumn{3}{l}{$r(x) =$ \tiny
$x^{24} + x^{23} - x^{21} - x^{20} + x^{18} + x^{17} - x^{15} - x^{14} + x^{12} - x^{10} - x^9 + x^7 + x^6
- x^4 - x^3 + x + 1$}\\
& \multicolumn{3}{l}{$t(x) = x^{14} - x + 1$}\\
& {$(a, b) = (0,12)$} & {$\dpol(X)=X^{26}+4$}
& {$x_0 = 2685463$}\\
\hline
\curvem
& \multicolumn{3}{l}{$p(x) = (x^{20} - x^{19} + x^{18} - x^{11} - 2x^{10} - x^9 + x^2 + 2x + 1)/3$}\\
& \multicolumn{3}{l}{$r(x) = (x^{18} - x^9 + 1)/3$}\\
& \multicolumn{3}{l}{$t(x) = -x^{10} + x + 1$} \\
& {$(a, b) = (0,9)$} & {$\dpol(X)=X^{27}+3$}
& {$x_0 = 374298113$}\\
\hline
\curvefv
& \multicolumn{3}{l}{$p(x) = 36x^4+36x^3+24x^2+6x+1$}\\
& \multicolumn{3}{l}{$r(x) = 36x^4+36x^3+18x^2+6x+1$}\\
& \multicolumn{3}{l}{$t(x) = 6x^2+1$}\\
& {$(a, b) = (0,7)$} & {$\dpol(X)=X^{12}+2$} & \\
& \multicolumn{3}{l}{$x_0 = 934494328215398161047821996449179050138683228531035586851830703422221130029$\textbackslash}\\
& \multicolumn{3}{l}{\hspace {6.7mm}$030240635045913079014$}
\end {tabular}
}
\caption {\label {table:curves256}
Curves for security level 256 bit}
\end {table}

\begin {table}
\centerline {
\begin{tabular}{c|c|l}
Curve      & Pairing & Vector\\
\hline
\curvez    & $\ev$   & $\left(\frac{1}{3}(x-2), \frac{1}{3}(x+1), \frac{1}{3}(x+1), -\frac{1}{3}(x+1), -\frac{1}{3}(x-2), -\frac{1}{3}(x-2)\right)$ \\
           & $\evt$  & $\left(\frac{1}{3}(x^3 + 1), \frac{1}{3}(2 - x^3)\right)$\\[0.5mm]
\hline               
\curveg    & $\ev$   & $\left(x^2, -x, 1, 0, 0, 0, 0, 0, 0, 0\right)$\\
\hline               
\curvefall & $\ev$   & $\left(6x + 2, 1, -1, 1\right)$ \\
           & $\evt$  & $\left(2x + 1, 6x^2 + 2x\right)$\\
\hline               
\curvep    & $\ev$   & $\left(1, 0, 0, 0, x, 0, 0, 0\right)$\\
           & $\evt$  & $\left(x^3 + x^2 - 1, x^4 + x^3 - x - 1\right)$\\
\hline               
\curves    & $\ev$   & $((2x-15)/35, -(11x-30)/35, -(2x-1)/7, (x+10)/35,$\\
           &         & $ (2x+5)/5,(8x+10)/35, (2x+6)/7, (17x+30)/35)), x=25\ \mathrm{mod}\ 70$\\
           & $\evt$  & $\left(49x^4 / 625, 7 + 168x^4 / 625\right)$\\
\hline               
\curvet    & $\ev$   & $\left(x^2, 1, 0, 0, 0, 0, 0, 0, 0, 0, 0, 0, 0, 0, 0, 0\right)$\\
\hline               
\curvei    & $\ev$   & $\left(1, 3x/7, 3x/7 + 1, 0, -2x/7, -2x/7-1\right)$\\
           & $\evt$  & $\left(18(x/7)^3 + 1, -(x/7)^3\right)$\\
\hline               
\curveh    & $\ev$   & $\left(x^2, -x, 1, 0, 0, 0, 0, 0, 0, 0, 0, 0, 0, 0, 0, 0, 0, 0, 0, 0\right)$\\
\hline               
\curveo    & $\ev$   & $\left(x^2, -x, 1, 0, 0, 0, 0, 0, 0, 0, 0, 0, 0, 0, 0, 0, 0, 0, 0, 0\right)$\\
\hline               
\curven    & $\ev$   & $\left(x^2, x, 1, 0, 0, 0, 0, 0, 0, 0, 0, 0\right)$\\
\hline               
\curvem    & $\ev$   & $\left(x, 0, 0, 0, 0, 0, 0, 0, 0, 0, 1, 0, 0, 0, 0, 0, 0, 0\right)$\\
           & $\evt$  & $\left(x^9, 1\right)$
\end{tabular}
}
\caption {\label {table:optivec}
Optimal pairings}
\end {table}
\begin {table}
\centerline {
\footnotesize
\begin{tabular}{l|c|c|p{9.6cm}}
  Curve       &   $s$ & $n$ & Addition sequence\\
  \hline
  \curvez     &  1    & 6   & [1, 2, 3]\\
  \hline
  \curveg     &   3   & 10  & [1, 2, 3, 4, 5, 6]\\
  \hline
  \curvefall  &   1   &  4  & [1, 2, \underline{3}, 6, 12, 18, 30, 36]\\
  \hline
  \hline
  \curvep     &  3    &  8 & [1, 2, 3, 4, 5, 6]\\
  \hline
  \curves     &  857500    &  8 & [\underline{1}, 2, 4, \underline{6}, 10, 11, 15, 20, 22, 25, 29, 30,
      40, 50, 54, 55, 75, 100, 125, 145, 220, 250, \underline{272}, 278, \underline{300}, \underline{440}, 585, 625,
      875, \underline{900}, \underline{950}, 1025, 1100, 1172, \underline{1226}, 1280, 1372,
      1390, 1750, \underline{1779}, 2290, 2780, 2925, 3000, 3300, \underline{4250},
      4375, 4704, 4750, 4850, 5125, \underline{9700}, 13000,
      \underline{13250}, 15000]\\
  \hline
  \curvet     &   4   & 16  & [1, 2, 3, 5]\\
  \hline
  \curvei     & 1029  &  6  & [\underline{1}, \underline{2}, 3, \underline{4}, 5, 7, 14, 15, 21, 25, 35, 49, 54, \underline{61}, 62, 70, 87, 98, 112, \underline{131},
    \underline{224}, 245, \underline{249}, 273, 319, 343, \underline{350}, \underline{364}, 434, 450, \underline{504}, 581, 609, 784, 931, \underline{1057},
    1407, \underline{1715}, 1911, \underline{2842}, \underline{4753}, 4802, 6517]\\
  \hline
  \curveh     &   3   & 18  & [1, 2, 3, 4]\\
  \hline
  \hline
  \curvel     & 3     &  8 & [1, 2, 3]\\
  \hline
  \curveo     & 3     & 20 & [1, 2, 3]\\
  \hline
  \curven     & 3     & 12 & [1, 2, 3, 4]\\
  \hline
  \curvem     & 3     & 18 & [1, 2, 3]
\end{tabular}
}
\caption {\label {table:addseq}
Addition sequences}
\end {table}

\clearpage
\paragraph {Acknowledgement.}
This research was partially funded by ERC Starting Grant ANTICS 278537.

\bibliographystyle{plain}
\bibliography{pairings}

\end{document}